\documentclass{article}

\usepackage{mytemplate}

\newtheorem{thm}{Theorem}
\newtheorem{prop}{Proposition}[section]

\newcommand{\U}{\mathbf{U}}
\newcommand{\SU}{{\mathbf{SU}}}
\renewcommand{\tr}{{\mathrm{tr}\,}}
\newcommand{\fq}{{\F_q}}
\newcommand{\fqr}{{\F_{q^r}}}
\newcommand{\fd}{{\mathcal{F}_d}}

\newcommand{\N}{\mathbf{N}}
\newcommand{\e}{\mathbf{e}}
\newcommand{\gt}{{\tilde{\gamma}}}

\begin{document}

\title{Artin-Schreier L-functions and Random Unitary Matrices}
\author{Alexei Entin\\ \\
Raymond and Beverly Sackler School of Math. Sciences\\
Tel Aviv University\\
Tel Aviv 69978\\
E-mail: alfxeyen@post.tau.ac.il}
\maketitle

\begin{abstract} We give a new derivation of an identity due to Z. Rudnick and P. Sarnak about the $n$-level correlations of eigenvalues of random unitary
matrices as well as a new proof of a formula due to M. Diaconis and P. Shahshahani expressing averages of trace products over the unitary matrix ensemble.
Our method uses the zero statistics of Artin-Schreier L-functions and a deep equidistribution result due to N. Katz.\end{abstract}

\section{Introduction}

In \cite{rudsar} Z. Rudnick and P. Sarnak computed the $n$-level correlation for the zeroes of the Riemann zeta function
(and in fact for a much larger class of L-functions) for a restricted class of test functions and obtained a complicated combinatorial expression involving the Fourier transform of the test function (see Theorem \ref{correlation} below). By a complicated and highly unstraightforward combinatorial argument they were able to show that this expression coincides with the $n$-correlation for the eigenvalues of random unitary matrices. 

In the present paper we will give a new derivation of this fact using a different approach. The main advantage of our new
approach is that it can be extended to cases to which the combinatorial approach has so far not been extended. The method has been applied
in \cite{hyper} to equate the $n$-level densities of quadratic L-functions and random unitary symplectic matrices for a class of test functions
unattained previously. Another derivation of the $n$-correlation for random unitary matrices has been recently given in \cite{conrey}.

Assume the Riemann Hypothesis. Denote the nontrivial zeroes of the Riemann zeta-function by $1/2+\gam_j,j=\pm 1,\pm 2,...$ so that 
$$\gam_j\in\R,\gam_{-j}=-\gam_j,|\gam_1|\le|\gam_2|\le....$$ Denote by $\gt_j=\gam_j\log|\gam_j|/2\pi$ the normalised zero ordinates, normalised so
that the mean spacing between the $\gt_j$ is 1. To define the $n$-correlation for the Riemann zeroes we fix a smooth test function $\phi:\R^n\to\C$ which is symmetric (i.e. unchanged by any permutation of the variables), satisfies $\phi(x_1+t,...,x_n+t)=\phi(x_1,...,x_n)$ for any $t\in\R$ and is a Schwartz function on the hyperplane
$\sum x_i=0$, i.e. $\phi$ and all its partial derivatives decay faster than any polynomial on this hyperplane.
We define the $n$-correlation of the first $T$ Riemann zeroes with test function $\phi$ to be
$$R_n(T,\phi)=\frac{1}{T}\sum_{1\le j_1,...,j_n\le T\atop{\mathrm{distinct}}}\phi(\gt_{j_1},...,\gt_{j_n}).$$
It is more convenient to consider the unrestricted sums
$$C_n(T,\phi)=\frac{1}{T}\sum_{1\le j_1,...,j_n\le T}\phi(\gt_{j_1},...,\gt_{j_n}).$$
The passage from unrestricted sums to restricted
ones is made by a standard technique called combinatorial sieving (essentially inclusion-exclusion), which expresses the restricted sum as a combination of the unrestricted sums for all $m\le n$
(with auxiliary test functions obtained from $\phi$ by identifying some of the variables). We omit the details as this passage is described in \cite[\S 4]{rudsar} among other places. We concentrate on the asymptotics of $C_n(T,\phi)$.

In \cite{rudsar} Z. Rudnick and P. Sarnak computed the limit of the $n$-correlation for the Riemann zeroes (and more general L-functions) as $T\to\ity$
for a restricted class of test functions. Suppose that we can write 
\beq\label{cond}\begin{split}\phi(x_1,...,x_n)=\int_{\R^n}\Phi(\xi_1,...,\xi_n)\del(\xi_1+...+\xi_n)e^{2\pi i\sum_{j=1}^n \xi_j x_j}\dd\xi_1...\dd\xi_n,\\ \sup\Phi\ss\{\textstyle\sum|\xi_j|<2\},\end{split}\eeq
(here $\del$ is the Dirac delta function) for some Schwartz function $\Phi\in\mathcal{S}(\R^n)$ which is supported on the set $\sum|\xi_j|<2$. 
Note that the values of $\Phi$ on the hyperplane $\xi_1+...+\xi_n=0$ are determined by $\phi$, it is essentially
the Fourier transform of $\phi$ restricted to $\sum x_j=0$. Our restriction on $\phi$ will be that $\Phi$ is supported on the set $\sum|\xi_j|<2$,
which we assume from now on.

We denote by $\e_i$ the standard basis vector $(0,...,1,...,0)\in\R^n$ (1 in the $i$-th position) and $\e_{i,k}=\e_i-\e_k$.
It is proved in \cite{rudsar} that under condition \rf{cond} (and for a fixed test function)
\begin{multline}\label{zetacor}\lim_{T\to\ity}C_n(T,\phi)=\\= \Phi(0)+\sum_{m=1}^\floor{n/2}\sum_{\{(\al_1,\be_1),...,(\al_m,\be_m)\}}\int_{\R^m}\Phi\lb\sum_{j=1}^n\xi_j\e_{\al_j,\be_j}\rb\prod_{j=1}^m |\xi_j|\dd\xi_j,\end{multline}
the sum being over all disjoint sets of pairs $\{(\al_1,\be_1),...,(\al_m,\be_m)\}$ with $1\le \al_j<\be_j\le n$.

Next we define the $n$-correlation for unitary matrices.
Let $N$ be a natural number, $\phi$ a test function as above. We define the associated periodic test function with scaling factor $N$ by
$$\tilde{\phi}(x_1,...,x_n)=\sum_{u_1,...,u_{n-1}\in\Z}\phi(N(x_1+u_1),...,N(x_{n-1}+u_{n-1}),Nx_n),$$ ($\Z$ denotes the set of integers). This function has period 1 in each variable.
Let $U$ be a size $N$ unitary matrix with eigenvalues $e^{2\pi i\th_j}, j=1,...,N$. The $\th_j$ are real numbers defined up to order and addition
of integers, so the quantity
$$R_n(U,\phi)=\frac{1}{N}\sum_{1\le j_1,...,j_n\le N\atop{\mathrm{distinct}}}\tilde{\phi}(\th_{j_1},...,\th_{j_n})$$
is well-defined, and we call it the $n$-correlation of $U$ with test function $\phi$.
Again it is simpler to consider the unrestricted sums
\beq\label{defc}C_n(U,\phi)=\frac{1}{N}\sum_{1\le j_1,...,j_n\le N}\tilde{\phi}(\th_{j_1},...,\th_{j_n}).\eeq 

By an ingenious combinatorial argument it was shown in \cite{rudsar} that $$\lim_{N\to\ity}\int_{\U(N)}C_n(U,\phi)\dd U$$ is also
given by the RHS of \rf{zetacor}, provided that the test function $\phi$ satisfies the condition \rf{cond}. Here the integration is w.r.t. the normalised Haar measure on $\U(N)$, the ensemble of $N\times N$ unitary matrices.
Consequently the $n$-correlations of zeta-zeroes and random unitary matrix eigenvalues coincide for this class of test functions
(it is conjectured that this is in fact true without any restriction on the test function)

Our goal is to give a new proof of the following result:
\begin{thm}\label{correlation} Let $\phi$ be a test function as above satisfying the condition \rf{cond}. Then
\begin{multline*}\int_{\U(N)}C_n(U,\phi)\dd U=\\=\Phi(0)+\sum_{m=1}^\floor{n/2}\sum_{\{(\al_1,\be_1),...,(\al_m,\be_m)\}}\int_{\R^m}\Phi\lb\sum_{j=1}^n\xi_j\e_{\al_j,\be_j}\rb\prod_{j=1}^m |\xi_j|\dd\xi_j+O(1/N)\end{multline*}
as $N\to\ity$, the sum being over all disjoint sets of pairs $\{(\al_1,\be_1),...,(\al_m,\be_m)\}$ with $1\le \al_j<\be_j\le n$. Here the implied constant in $O(1/N)$ may
depend on $\phi$.
\end{thm}
The main term appearing here is exactly the RHS of \rf{zetacor}.

A closely related quantity to the $n$-level correlations of the eigenvalues of random unitary matrices is the average trace product, namely
\beq\label{defm}M(r_1,...,r_n;N)=\int_{\U(N)}\prod_{i=1}^n\tr U^{r_i}\dd U,\eeq where $r_1,...,r_n$ are integers, $\U(N)$ denotes the ensemble of $N\times N$ unitary
matrices and integration is w.r.t. the normalised Haar measure on $\U(N)$. The following theorem is proved by P. Diaconis and M. Shahshahani in \cite{diac}. Our statement is taken from \cite{diac2}.

\begin{thm}\label{main} Let $r_1,...,r_n,N$ be nonzero integers s.t. $\sum|r_i|\le 2N$. Let $s_1,...,s_m$ be the distinct values appearing in the list 
$|r_i|,i=1,...,n$ and $a_j,b_j,j=1,...,m$ be such that $s_j$ appears $a_j$ times in the list $r_1,...,r_n$ while $-s_j$ appears $b_j$ times.
Then $$M(r_1,...,r_n;N)=\choice{\prod_{j=1}^m a_j!s_j^{a_j},}{a_j=b_j, j=1,...,m,}{0,}{\mbox{otherwise.}}$$\end{thm}

We will derive Theorem \ref{correlation} from Theorem \ref{main} by a standard calculation which will be carried out in section \ref{seccor}.
Our method for proving Theorem \ref{main} proceeds by computing the average of $\prod_{i=1}^n\tr U^{r_i}$ over the set of Frobenius elements of a suitable family of Artin-Schreier
L-functions and using the equidistribution result of N. Katz \cite[Theorem 3.9.2]{katz} to relate this to the unitary matrix average. A similar approach
was applied in \cite{hyper} to the unitary symplectic ensemble using hyperelliptic curves. The method of \cite{diac} is completely different
and is based on the representation theory of unitary matrices.

\section{Artin-Schreier L-functions}

In this section we summarise the properties of Artin-Schreier (A-S in short) L-functions that we will need. A more detailed survey and further
references may be found in \cite[\S 3,\S 7]{as}. An exposition of the basic properties of A-S L-functions with full details may be found in \cite{stepanov} and of
Dirichlet L-functions for the ring of polynomials (which will be needed shortly) in \cite{rosen}.

Let $p$ be a prime number $q$ its power, $\fq$ the finite field with $q$ elements. Let $f\in\fq[x]$ be a polynomial of degree $d$ prime to $p$.
For the rest of the paper we always assume $(d,p)=1$.
Fix an additive character $\psi:\F_p^+\to\C^\times$.
The Artin-Schreier (A-S) function with defining polynomial $f$ is given by 
\beq\label{lfunc}L_f(z)=\exp\lb\sum_{r=1}^\ity\sum_{\al\in\F_{q^r}}\psi\lb\mathrm{tr}_{\fqr/\F_p}f(\al)\rb\frac{z^r}{r}\rb.\eeq
It is known that $L_f(z)$ is in fact a polynomial of degree $d-1$ and in fact
$$L_f(z)=\prod_{j=1}^{d-1}(1-q^{1/2}e^{2\pi i\th_j}z),$$ where $\th_j=\th_{f,j}$ are real numbers (well defined upto the addition of integers and reordering). This is an equivalent formulation of the Riemann Hypothesis
for the curve $y^p-y=f(x)$ over $\fq$. We denote by $\Th_f$ the conjugacy class of the matrix 
$\diag(e^{2\pi i\th_{f,1}},...,e^{2\pi i\th_{f,{d-1}}})\in\U(d-1)$. It is called the \emph{Frobenius class} of $\Th_f$.

Denote $$\fd=\{f\in\F_q[x]|\deg f=d, f(0)=0\}.$$ The family of L-functions $\{L_f\}_{f\in\fd}$ is independent of the choice of additive character
$\psi$, because replacing $\psi$ with $\psi^a$ (with $(a,p)=1$) has the same effect on the $L$-function as replacing $f$ with $af$ (by \rf{lfunc}).
The following deep result due to N. Katz (special case of \cite[Theorem 3.9.2]{katz}) is the main ingredient of the present work:

\begin{thm}\label{equid} Let $d$ be fixed and assume $p>5$. Then as $q\to\ity$ the family $\{\Th_f\}_{f\in\fd}$ becomes equidistributed in the space of conjugacy classes
of some Lie group $\SU(d-1)\ss G\ss\U(d-1)$ with the measure induced from the Haar measure on $G$.\end{thm}

We will also make use of the connection between A-S L-functions and Dirichlet L-functions. If $Q(x)\in\F_q[x]$ is a monic polynomial and
$\chi$ is a Dirichlet character modulo $Q$ we define its L-function 
$$L_\chi(z)=\sum_{u\in\fq[x]\atop\mathrm{monic}}\chi(u)z^{\deg u}=\prod_{P\in\fq[x]\atop{\mathrm{prime}}}\lb 1-\chi(P)z^{\deg P}\rb^{-1}.$$
A Dirichlet character $\chi$ is called even if it is trivial on $\F_q$. It is called primitive if it is not induced from a character of
smaller modulus.

A basic fact we will use is the following (see \cite[\S 7]{as} for a proof):

\begin{prop} To each $f\in\fd$ we can assign a primitive even Dirichlet character $\chi_f$ modulo $Q=x^{d+1}$ s.t. $L_f(z)=(1-z)^{-1}L_{\chi_f}(z)$.
If $p>d$ then this gives a bijection between $\fd$ and the set of even primitive Dirichlet characters $\chi$ modulo $Q$.\end{prop}

We will also make use of the \emph{explicit formula} (a proof can be found in \cite[\S 7]{as} as well):
\beq\label{expform}\tr\Th_f^r=-q^{-r/2}-q^{-r/2}\sum_{u\atop{\mathrm{monic}\atop{\deg u=r}}}\Lam(u)\chi_f(u).\eeq
Here $r$ is any natural number and $$\Lam(u)=\choice{\deg P,}{u=P^k\mbox{ for some prime } P,}{0,}{\mbox{otherwise}}$$
is the von Mangoldt function.

A final fact that we will need is the orthogonality relation for characters. Let $Q=x^{d+1}$ be our modulus, $u\in\fq[x]$ nonconstant and prime to $x$.
Then
\beq\label{ort}\av{\chi(u)}_\chi=\left\{\begin{array}{ll}1, & u\equiv a\pmod{x^{d+1}}\mbox{ for some }a\in\fq^\times,\\
-1/(q-1), & u\equiv a+bx^d\pmod{x^{d+1}}\mbox{ for some }a,b\in\fq^\times,\\ 0, & \mbox{otherwise,}\end{array}\right.\eeq
where $\av{\cdot}_\chi$ denotes the average taken over all even primitive characters modulo $Q$.

\section{Average trace products for Artin-Schreier L-functions}\label{asav}

We keep the notation of the previous section. In particular $q$ is a power of a prime $p$ and $d$ is a natural number prime to $p$. We also
assume $p>d$ for this section. In the present section we will compute an estimate for the quantity
$$\av{\prod_{i=1}^k\tr\Th_f^{r_j}\prod_{j=1}^l\tr\Th_f^{-t_j}}_{f\in\fd},$$ where $r_1,...,r_k,t_1,...,t_l$ are natural numbers satisfying 
$\sum r_i=\sum t_j <d$. We will see that the case $\sum r_i=\sum t_j$ is really all we need. In the following
section we will combine this estimate with Theorem \ref{equid} to obtain Theorem \ref{main}.

For the rest of this section the asymptotic big-$O$ notation will always have an implicit constant which may depend on $k,l,d$ (which we assume are fixed),
but not on $p,q$. We begin by applying \rf{expform}. Take some $f\in\fd$. We have

\begin{multline*}\prod_{i=1}^k\tr\Th_f^{r_i}\prod_{j=1}^l\tr\Th_f^{-t_j}=\prod_{i=1}^k\tr\Th_f^{r_i}\prod_{j=1}^l\overline{\tr\Th_f^{t_j}}=\\
=(-1)^{k+l}q^{-\sum r_i}\sum_{u_1,...,u_k,v_1,...,v_l\atop{\mathrm{monic}\atop{\deg u_i=r_i,\deg v_j=t_j}}}
\prod_{i=1}^k\Lam(u_i)\chi_f(u_i)\prod_{j=1}^l\Lam(v_j)\overline{\chi}_f(v_j) + O(q^{-1/2})\end{multline*}
(the error term comes from the term $q^{-r/2}$ in \rf{expform} and the fact that $\tr\Th_f^r=O(1)$).

Now using the orthogonality relation \rf{ort} and the fact that for $p>d$ the characters $\chi_f$ for $f\in\fd$ are precisely all the even primitive characters
modulo $x^{d+1}$ we conclude that:

\begin{multline*}\av{\prod_{i=1}^k\tr\Th_f^{r_i}\prod_{j=1}^l\tr\Th_f^{t_j}}_{f\in\fd}=\\=(-1)^{k+l}q^{-\sum r_i}\sum_{a\in\F_q^\times}\sum_{u_1,...,u_k,v_1,...,v_l
\atop{\mathrm{monic}\atop{(u_iv_j,x)=1\atop{\deg u_i=r_i,\deg v_j=t_j\atop{u_1...u_kv_1^{-1}...v_l^{-1}\equiv a\pmod{x^{d+1}}}}}}}\prod_{i=1}^k\Lam(u_i)\prod_{j=1}^l\Lam(v_j)-\\-(-1)^{k+l}\frac{q^{-\sum r_i}}{q-1}
\sum_{a,b\in\F_q^\times}\sum_{u_1,...,u_k,v_1,...,v_l
\atop{\mathrm{monic}\atop{(u_iv_j,x)=1\atop{\deg u_i=r_i,\deg v_j=t_j\atop{u_1...u_kv_1^{-1}...v_l^{-1}\equiv a+bx^d\pmod{x^{d+1}}}}}}}\prod_{i=1}^k\Lam(u_i)\prod_{j=1}^l\Lam(v_j)+\\+O(q^{-1/2}).\end{multline*}

Now since $\sum r_i=\sum t_j<d$, for $u_i,v_j$ monic and prime to $x$ with $\deg u_i=r_i,\deg v_j=t_j$ we can only have $$\prod u_i\equiv (a+bx^d)\prod v_j\pmod{x^{d+1}}$$
for $a\in\F_q^\times,b\in\F_q$ if $a=1,b=0$ and $\prod u_i=\prod v_j$. Therefore we have
\begin{multline}\label{main1}\av{\prod_{i=1}^k\tr\Th_f^{r_i}\prod_{j=1}^l\tr\Th_f^{-t_j}}_{f\in\fd}=\\=(-1)^{k+l}q^{-\sum r_i}\sum_{u_1,...,u_k,v_1,...,v_l
\atop{\mathrm{monic}\atop{(u_iv_j,x)=1\atop{\deg u_i=r_i,\deg v_j=t_j\atop{u_1...u_k=v_1...v_l}}}}}\prod_{i=1}^k\Lam(u_i)\prod_{j=1}^k\Lam(v_j)+O(q^{-1/2}).\end{multline}

Now the contribution to \rf{main1} of $u_i,v_j$ some of which are proper prime powers or such that two of the $u_i$ or two of the $v_j$ coincide is easily seen to be $O(q^{-1})$.
For example suppose that $u_1$ is a power (higher than 1) of some polynomial. Then there are at most $q^{r_1-1}$ possibilities for $u_1$ and
at most $q^{r_i}$ possibilities for every other $u_i$. Once the $u_i$ are determined there are at most $d^n=O(1)$ ways to factor the product $\prod u_i$
into $l$ factors $v_j$ (because each $v_j$ is composed of a subset of the prime factors of $\prod u_i$ of which there are at most $d$).
After multiplying by $q^{-{\sum r_i}}$ and using $\Lam(u_i),\Lam(v_j)\le d$ we get a total contribution of at most $O(q^{-1})$.
Coincidences of the form $u_1=u_2$ also contribute at most $O(q^{-1})$ and this is seen similarly.

The main contribution is thus from sets of primes $u_1,...,u_k,v_1,...,v_l$ s.t. the $u_i$ are distinct and the $v_j$ are distinct. When such $u_i$ are chosen the $v_j$ are the same as the $u_i$ up to a change of order (in particular we get this contribution only if $k=l$ and the $r_i$ and $t_j$ coincide up to order). The number
of tuples of distinct primes $u_1,...,u_k$ with $\deg u_i=r_i$ is $q^{\sum r_i}/\prod{r_i}+O\lb q^{\sum r_i-1}\rb$ (this follows from the fact that the number of prime polynomials of degree $r$ is $q^r/r+O\lb q^{\floor{r/2}}/r\rb$). Now as in the statement of Theorem \ref{main} let $s_1,...,s_m$ be the distinct
values appearing in the list $r_1,...,r_k$, each appearing $a_i$ times. If $k=l$ and the $r_i=t_i$ then for each choice of $u_i$ we have $\prod_{\nu=1}^m a_\nu!$ ways to order the
$v_j$ (so that $\deg v_j=t_j$), so taking everything together we
conclude that \beq\label{est}\av{\prod_{i=1}^k|\tr\Th_f^{r_i}|^2}_{f\in\fd}=\prod_{j=1}^m a_j!\prod_{i=1}^k t_i+O(q^{-1/2}).\eeq
On the other hand if the $r_i$ and the $t_j$ do not coincide up to order (e.g. if $k\neq l$) then
\beq\label{est2}\av{\prod_{i=1}^k\tr\Th_f^{r_i}\prod_{j=1}^l\tr\Th_f^{-t_j}}_{f\in\fd}=O(q^{-1/2}).\eeq

\section{Proof of Theorem \ref{main}}

We are now ready to prove Theorem \ref{main}. We keep the notation of the previous two sections and assume $p>d$. We fix natural numbers $n,d,r_1,...,r_k$ with $\sum r_i<d$ and take 
the limit $q\to\ity$. Then by Theorem \ref{equid} we have
$$\av{\prod_{i=1}^n|\tr\Th_f^{r_i}|^2}_{f\in\fd}\to\int_{\U(d-1)}\prod_{i=1}^n|\tr U^{r_i}|^2=M(r_1,...,r_k,-r_1,...,-r_k;d-1)$$ as $q\to\ity$
(using the notation of section 1 and the integral being as usual w.r.t. the normalised Haar measure). We may integrate over $\U(d-1)$ and not some smaller Lie group
$\SU(d-1)\ss G\ss\U(d-1)$ as stated in Theorem \ref{equid} because unitary scalars make no difference to the absolute value of the traces.
Using the estimate \rf{est} we obtain
$$M(r_1,...,r_k,-r_1,...,-r_k;d-1)=\prod_{j=1}^m a_j!\prod_{i=1}^k r_i,$$ where the $s_j,a_j$ are determined from $r_i$ as in the end of the previous
section. This settles Theorem \ref{main} for $r_i$ that can be ordered so that $r_{k+i}=-r_{i},1\le i\le k$ and $n=2k$ is even, where we take $d=N-1$.
Similarly using \rf{est2} and Theorem \ref{equid} we obtain (for natural numbers $r_i,t_j$).
$M(r_1,...,r_k,-t_1,...,-t_l;N)=0$ if $\sum r_i=\sum t_j$ but the $t_j$ are not a reordering of the $r_i$.

The only remaining case to consider is $M(r_1,...,r_k,-t_1,...,-t_l;N)$ (again the $r_i,t_j$ are natural numbers) where $\sum r_i\neq\sum t_j$.
But the map $U\mapsto e^{i\al}U$ (with $\al\in\R$) preserves the Haar measure on $\U(N)$ and multiplies the trace product $\prod_{i=1}^k\tr U^{r_i}\prod_{j=1}^l\tr U^{-t_j}$ by
$e^{i\al(\sum r_i-\sum t_j)}$ and so $$M(r_1,...r_k,-t_1,...,-t_l;N)=e^{i\al(\sum r_i-\sum t_j)}M(r_1,...r_k,-t_1,...,-t_l;N)$$
for every $\al\in\R$. It follows that $M(r_1,...r_k,-t_1,...,-t_l;N)=0$ if $\sum r_i\neq\sum t_j$.
Note that the condition $\sum r_i+\sum t_j\le 2N$ is not required in this case.

\section{$n$-correlations} \label{seccor}

In the present section we derive Theorem \ref{correlation} from Theorem \ref{main}. We do this by a standard calculation involving Fourier series.
For the rest of the section we fix a natural number $n$ and a test function $\phi$ satisfying the assumptions made in the introduction. 
That is $\phi:\R^n\to\ity$ is symmetric (unchanged by any permutation of the variables), translation invariant in the sense that
$\phi(x_1+t,...,x_n+t)=\phi(x_1,...,x_n)$ for all $t\in\R$ and can be expressed as
$$\phi(x_1,...,x_n)=\int_{\R^n}\Phi(\xi_1,...,\xi_n)\del(\xi_1+...+\xi_n)e^{2\pi i\sum_{j=1}^n \xi_i x_i}\dd\xi_1...\dd\xi_n$$
($\del$ denotes the Dirac delta function, the values of $\Phi$ on the hyperplane $\sum\xi_j=0$ are determined by $\phi$) where
$\Phi$ is a Schwartz function supported on $\sum|\xi_j|<2$.
We also define the periodic test function associated with $\phi$ with scaling factor $N$ by $$\tilde{\phi}(x_1,...,x_n)=\sum_{u_1,...,u_{n-1}\in\Z}\phi(N(x_1+u_1),...,N(x_{n-1}+u_{n-1}),Nx_n).$$

The Fourier series expansion of $\tilde{\phi}$ is the following:
$$\tilde{\phi}(x_1,...,x_n)=\frac{1}{N^{n-1}}\sum_{r_1,...,r_n\in\Z\atop\sum r_j=0}\Phi\lb\frac{r_1}{N},...,\frac{r_n}{N}\rb e^{2\pi i\sum_{j=1}^nr_j x_j}.$$
Combined with \rf{defc} we get
$$C_n(U,\phi)=\frac{1}{N^n}\sum_{r_1,...,r_n\in\Z\atop{\sum r_j=0}}\Phi\lb\frac{r_1}{N},...,\frac{r_n}{N}\rb\prod_{j=1}^n\tr(U^{r_j}).$$
Averaging and using \rf{defm} we get
\beq\label{avc}\int_{\U(N)}C_n(U,\phi)\dd U=\frac{1}{N^n}\sum_{r_1,...,r_n\in\Z}\hat{\phi}\lb\frac{r_1}{N},...,\frac{r_n}{N}\rb M(r_1,...,r_n;N)\eeq
(we omitted the condition $\sum r_j=0$ because otherwise we have $M(r_1,...,r_n;N)=0$ by Theorem \ref{main}).

At this point we use the assumption that $\Phi(\xi_1,...,\xi_n)$ is supported on the set $\sum_{i=1}^n|\xi_i|<2$. Thus the sum in \rf{avc}
is only over $r_1,...,r_n$ satisfying $\sum|r_j|<2N$ (and $\sum r_j=0$). As in the statement of Theorem \ref{main}, for each
tuple $r_1,...,r_n$ of integers we consider the distinct nonzero values $s_1,...,s_m$ appearing among the $r_i$ (note that some $r_i$ may be zero).
If $s_j$ appears $a_j$ times among the $r_i$ and $-s_j$ appears $b_j$ times then $\sum a_j+\sum b_j\le n$ and nonzero contribution comes only from
the case $a_j=b_j$ in which case 
\beq\label{m}M(r_1,...,r_n;N)=N^{n-2\sum a_j}\prod_{j=1}^n a_j!s_j^{a_j}\eeq
(we used the fact that $\tr U^0=N$).

We make the convention that the asymptotic big-$O$ notation has an implicit constant which may depend on $n$ and $\phi$ but not on $N$ (note the difference
from the convention adopted in the previous section, where $N$ was also assumed fixed!). First we observe that the case in which $a_j=b_j>1$ for some $j$ contributes only $O(1/N)$
to \rf{avc}. Indeed for any choice of $a_1,...,a_m$ (and $b_j=a_j$) we have $O(N^m)$ choices for the $s_j$ (satisfying $s_j\le N$) and each contributes
$O\lb N^{-\sum a_j}\rb$ (by \rf{m} and due to the factor $N^{-n}$ appearing in \rf{avc}), so we only get a significant contribution from
those tuples with $\sum a_j=1$, i.e. $a_j=b_j=1$, while the rest contribute $O(1/N)$.

Now let $m\le\floor{n/2}$ be a natural number and let $\sig=\{(\al_1,\be_1),...,(\al_m,\be_m)\}$ be a set of disjoint distinct pairs of indices 
$1\le \al_j,\be_j\le n$.
For each such $\sig$ consider the contribution to \rf{avc} of the tuples $r_1,...,r_n$ of integers s.t. $r_{\al_j}=-r_{\be_j}>0$ and $r_i=0$ if $i$ is not among
the $\al_j,\be_j$. We denote by $\e_i$ the standard basis vector $(0,...,1,...,0)\in\R^n$ (1 in the $i$-th position) and $\e_{i,k}=\e_i-\e_k$.
The contribution to \rf{avc} from the considered tuples is (by \rf{m})

$$\frac{1}{N^{2m}}\sum_{s_1,...,s_m\in\N}\Phi\lb\frac{1}{N}\sum_{j=1}^ns_j\e_{\al_j,\be_j}\rb\prod_{j=1}^m s_j$$
($\N$ denotes the set of natural numbers). This is a Riemann sum (with step $1/N$) approximating up to $O(1/N)$ the integral
$$\int_{\R_+^m}\Phi\lb\sum_{j=1}^n\xi_j\e_{\al_j,\be_j}\rb\prod_{j=1}^m \xi_j\dd\xi_j.$$
We conclude that 

\begin{multline*}\av{C_n(U,\phi)}_{U\in\U(N)}=\\=\Phi(0)+\sum_{m=1}^\floor{n/2}\sum_{\{(\al_1,\be_1),...,(\al_m,\be_m)\}}\int_{\R_+^m}\Phi\lb\sum_{j=1}^n\xi_j\e_{\al_j,\be_j}\rb\prod_{j=1}^m \xi_j\dd\xi_j+O(1/N),\end{multline*}
where the sum is over all the disjoint sets of distinct pairs $\{(\al_1,\be_1),...,(\al_m,\be_m)\}$. The term $\Phi(0)$ is the contribution of
$r_1=...=r_n=0$ to \rf{avc}.

If we use the symmetry (invariance under permutation of variables) of $\phi$ (implying the symmetry of $\Phi$ restricted to $\sum\xi_i=0$) we
can rewrite this as
\begin{multline*}\av{C_n(U,\phi)}_{U\in\U(N)}=\\=\Phi(0)+\sum_{m=1}^\floor{n/2}\sum_{\{(\al_1,\be_1),...,(\al_m,\be_m)\}}\int_{\R^m}\Phi\lb\sum_{j=1}^n\xi_j\e_{\al_j,\be_j}\rb\prod_{j=1}^m |\xi_j|\dd\xi_j+O(1/N),\end{multline*}
the sum now being over disjoint sets $\{(\al_1,\be_1),...,(\al_m,\be_m)\}$ with $\al_j<\be_j$. This concludes the proof of Theorem \ref{correlation}.

{\bf Acknowledgment.} The author would like to thank Ze\'{e}v Rudnick for suggesting this line of research and for many valuable discussions and suggestions in the course of research and writing the paper. The present work is part of the author's Ph.D. studies in Tel-Aviv University under his supervision.
The author would also like to thank Chantal David for pointing out some minor errors in previous versions of the paper.

The research leading to these results has received funding from the European
Research Council under the European Union's Seventh Framework Programme
(FP7/2007-2013) / ERC grant agreement no. 320755.

\end{document}